%% file: expgrowthNY.tex
\documentclass[11pt]{article}
\input{macros_angl_pdf.tex}

\usepackage{mathtools}

\title{From exponential counting to pair correlations\\
}
\author{Jouni Parkkonen \and Fr\'ed\'eric Paulin} 
\date{\today}

\begin{document}
\bibliographystyle{../alphanum}	
\maketitle

\begin{abstract}
We prove an abstract result on the correlations of pairs of elements
in an exponentially growing discrete subset $\E$ of $[0,+\infty[\,$
endowed with a weight function.  Assume that there exist
$\alpha\in\RR$, $c,\delta>0$ such that, as $t\ra+\infty$, the weighted
number $\wt\omega(t)$ of elements of $\E$ that are not greater than
$t$ is equivalent to $c\,t^\alpha e^{\delta t}$. We prove that the
distribution function of the unscaled differences of elements of $\E$
is $t\mapsto\frac\delta 2\,e^{-|t|}$, and that, under an error term
assumption on $\wt\omega(t)$, the pair correlation with a scaling with
polynomial growth exhibits a Poissonian behaviour.  We apply this
result to answer a question of Pollicott and Sharp on the pair
correlations of closed geodesics and common perpendiculars in
negatively curved manifolds and metric graphs. 
\end{abstract}

\section{Introduction}

When studying the asymptotic distribution of a sequence of finite
subsets of $\RR$, finer information is sometimes given by the
statistics of the spacing between pairs or $k$-tuples of elements,
seen at an appropriate scaling. This problematic is largely developped
in quantum chaos, including energy level spacings or clusterings, and
in statistical physics, including molecular repulsion or interstitial
distribution.  See for instance \cite{Montgomery73, Berry88, RudSar98,
BocZah05, MarStr13, LarSto20, HofKal21}. In \cite{PolSha06,PolSha13},
Pollicott-Sharp study the pair correlations of lengths of closed
geodesics in negatively curved manifolds as the word length of the
elements of the fundamental group that represent them tends to
$+\infty$. They mention that a result replacing the word length by the
Riemannian length does not seem to be available. One aim of this note
is to answer this problem, by a very general method.

For any set $\E$, a {\em weight function} (or multiplicity function
when its values are positive integers) on $\E$ is simply a function
$\omega:\E\ra\; ]\,0,+\infty\,[\,$.  The {\em growth function} (or
counting function when the weights are integers) of a locally finite subset
$\E$ of $[\,0,+\infty\,[$ endowed with a weight function $\omega$ is
the map $\N_{\E,\omega}:[0,+\infty\,[\ra\; [\,0,+\infty\,[\,$ defined
by
\[
\N_{\E,\omega}:t\mapsto \sum_{x\in\E\cap[0,t]} \omega(x)\;.
\]
Let $\F=(\,(F_N)_{N\in\NN}, \;\omega)$ be a nondecreasing sequence of
finite subsets $F_N$ of a finite dimensional Euclidean space $E$,
endowed with a weight function $\omega: \bigcup_{N\in\NN} F_N \ra\;
]\,0,+\infty\,[\,$. Let $\psi$ be any function from $\NN$ to
$[1,+\infty\,[$, called a {\em scaling function}, and let $\psi':
\NN\ra]0,+\infty[$ be an appropriately chosen function, called a
{\em renormalising function}. The {\em pair correlation measure
of $\F$ at time $N$ with scaling $\psi(N)$} is the measure on $E$ with
finite support
\begin{equation}\label{eq:defiPCmeasure}
\R^{\F,\psi}_N=\sum_{x,y\in F_N}\;\omega(x)\,\omega(y)\,
\Delta_{\psi(N)(y-x)}\;,
\end{equation}
where $\Delta_z$ denotes the unit Dirac mass at $z$ in any measurable
space.  When the sequence of measures $\R^{\F,\psi}_N$, renormalised
by $\psi'(N)$, weak-star converges to a measure $g_{\F,\psi}\,\Leb_E$
absolutely continuous with respect to the Lebesgue measure $\Leb_E$ of
$E$, the Radon-Nikodym derivative $g_{\F,\psi}$ is called the
asymptotic {\it pair correlation function} of $\F$ {\em for the
  scaling} $\psi$ {\em and renormalisation} $\psi'$. When
$g_{\F,\psi}$ is a positive constant, we say that $\F$ has a {\em
  Poissonian behaviour for the scaling} $\psi$ {\em and
  renormalisation} $\psi'$.

\btheo\label{theo:mainintro} Let $\E$ be a locally finite subset of
$[0,+\infty[\,$ endowed with a weight function $\omega$.  Assume that
there exist $\alpha\in\RR$, $c,\delta>0$ and $\kappa\geq 0$ such
that, as $t\ra+\infty$, we have
\[
\N_{\E,\,\omega}(t)\sim c\;t^\alpha\,e^{\delta \,t}(1+o(e^{-\kappa\, t}))
\;.
\]
Let $\psi:\NN\ra[1+\infty[$ be an at most polynomially growing scaling
function, with renormalising function $\psi': N\mapsto
\frac{\N_{\E,\,\omega}(N)^2}{\psi(N)}$. Then the family $\F=
(\,(F_N=\{x\in \E:x\leq N\})_{N\in\NN}, \;\omega)$ has a pair
correlation function $g_{\F,1}:t\mapsto\frac{\delta}{2}\;
e^{-\,\delta\,|t|}$ if $\psi= 1$, and has Poissonian behaviour with
$g_{\F,\psi}= \frac{\delta}{2}$ if ${\displaystyle\lim_{+\infty}\psi
  =\infty}$ and $\kappa>0$.  \etheo

\medskip
We give some comments on the above statement at the beginning of
Section \ref{sec:proof}. We refer to Theorem \ref{theo:main} for a
more precise version, including error terms. The work on error terms
constitutes the main technical parts of this paper.

Numerous settings in number theory, in geometry and in dynamical
systems\footnote{See, for instance \cite{ParPol83}, \cite{EskMir11}.}
give rise to counting functions that satisfy the assumption of Theorem
\ref{theo:mainintro}. We will give some applications of the above
result on geometry and dynamics in Section \ref{sec:application}.
Following the notation of \cite{PolSha06}, for all $a<b$ in $\RR$ and
$N\in\NN$, let
\[
\pi_\E(N,[a,b])=\R^{\F,1}_N([a,b])=
\sum_{x,y\in \E\;:\;x,y\leq N,\;a\leq y-x\leq b}\;\omega(x)\,\omega(y)
\]
be the weighted number of differences of elements in $\E\cap [0,N]$
that lie in the interval $[a,b]$. Since the limit measure is atomless,
under the assumptions of Theorem \ref{theo:main}, we have the
following corollary (see also Corollary \ref{cor:geometric}).

\bcoro
For all $a<b$ in $\RR$, as $N\ra+\infty$, we have
$$
\pi_\E(N,[a,b])\sim \frac\delta2\;\N_{\E,\,\omega}(N)^2
\int_a^be^{-\delta|t|}\;dt\;.\;\;\;\Box
$$
\ecoro

This answers the question of Pollicott-Sharp \cite{PolSha06} when $\E$
is the set of lengths of closed geodesics in a closed negatively
curved manifold, $\delta$ is the topological entropy of its geodesic
flow and $\omega$ is the multiplicity function of these lengths (see
Remark (3) in Section \ref{sec:proof}).

We suspect that when the scaling function has superexponential growth,
the empirical measures $\frac{1}{\psi'(N)}\R^{\F,\psi}_N$ have a total
loss of mass as $N\ra+\infty$ whatever the renormalising function
$\psi'$ is, hence that the pair correlation function $g_{\F,\psi}$
exists and is identically $0$. The main open problem related to
Theorem \ref{theo:mainintro} is to study the pair correlations for
scaling functions $\psi$ which are at the threshold, that is, are just
exponentially growing. For instance, the set $\E= \{\,\ln n :
n\in\NN-\{0\}\}$ endowed with the trivial multiplicity function
$\omega:x\mapsto 1$ satisfies the assumption of the above theorem with
$c=1$, $\alpha=0$ and $\delta=1$. In \cite{ParPau22a}, we study the
pair correlations of this family for general scalings and some
arithmetic weights functions, proving surprising level repulsion
phenomena when $\psi(N)=e^N$.

\medskip\noindent
{\small {\it Acknowledgements: } This research was supported by the
  French-Finnish CNRS IEA BARP. }

\section{Preliminaries on the growth of positive sequences}
\label{sect:prelimposseq}

In this section, we recall some standard terminology used in the
paper, and we prove two technical results used in the proof of the
main results in Section \ref{sec:proof}.
  
Recall that given a set of parameters $P$ and a positive map $h$
defined on a neighborhood of $+\infty$ in $\NN$ or $\RR$, we denote by
$\bigO_P(h)$ (respectively $\smallo_P$) any {\it Landau function} (as
the variable goes to $+\infty$) from $\RR$ to $\RR$ such that there
exists a constant $M>0$ depending only on the parameters in $P$ and
$t_0\geq 0$ (possibly depending on ambient data) such that for every
$t\geq t_0$, we have $|\bigO_P(h)(t)| \leq M\;h(t)$ (respectively such
that ${\displaystyle\lim_{+\infty} }\;
\frac{|\smallo_P(h)(t)|}{h(t)}=0$).

A positive sequence $(x_n)_{n\in\NN}$ is

$\bullet$~ {\it subexponentially growing} if for every
$\ga>0$, we have $\lim_{n\ra+\infty} \frac{x_n}{e^{\ga n}}=0$,

$\bullet$~ {\it at most polynomially growing} if there exists $\ga>0$
such that $\lim_{n\ra+\infty} \frac{x_n}{n^{\ga}}=0$,

$\bullet$~ {\it strictly sublinearly growing} if there exists
$\ga\in\;]0,1[$ such that $\lim_{n\ra+\infty} \frac{x_n}{n^{\ga}}=0$.

\medskip
The first result generalises a classical result on the geometric sums
(when $b=0$) to the generality needed for the proofs in Section
\ref{sec:proof}.

\blemm\label{lem:gengeomlimit}
For every $b\in\RR$, for every sequence $(a_M)_{M\in\NN}$ in
$]1,+\infty[$ such that the sequence $(\frac{1} {\ln a_M} )_{M\in\NN}$
is strictly sublinearly growing, as $M$ tends to $+\infty$ in
$\NN$, we have
$$
\sum_{k=1}^M k^b (a_M)^{k}=
\frac{a_M}{a_M-1}M^{b} (a_M)^{M}\big(1+\bigO_b(\frac{1}{\sqrt{M}})\big)\;.
$$
\elemm

\dem Let $\ga\in\;]0,1[$ be such that $\lim_{M\ra+\infty} M^{\ga}\ln
a_M=+\infty$. As a preliminary remark, note that we have
$n^b (a_M)^{n} = \bigO_b(M^b (a_M)^{M})$ for every $n\in \{1,\dots,
M\}$~: This is immediate if $b \geq 0$, and follows when $b < 0$ by
considering separately the case $n\geq \frac{M}{2}$ (in which case we
have $\frac{n^b(a_M)^{n}}{M^b (a_M)^{M}} \leq 2^{|b|}$) and $n\leq
\frac{M}{2}$ (in which case we have
$$
\frac{n^b (a_M)^{n}}{M^b (a_M)^{M}} \leq M^{|b|}(a_M)^{-M/2}=
e^{-\frac{M^{1-\ga}}{2}\big(M^\ga\ln a_M-2|b|\frac{\ln M}{M^{1-\ga}}\big)}\;,
$$
which converges to $0$ as $M$ tends to $+\infty$).

Recall that $(1-\frac{1}{n+1})^b =1+\bigO_b(\frac{1}{n})$ for every
$n\geq 1$. With $\Sigma_M= \sum_{n=1}^M n^b (a_M)^{n}$, for every
$S\in[1,M]$, by the standard telescopic sum argument and by cutting
the third sum below for $n\leq S$ and for $n> S$ (using in this second
case the preliminary remark), we hence have
\begin{align*}
&(a_M-1)\;\Sigma_M= \sum_{n=1}^M n^b  (a_M)^{n+1} -\Sigma_M=
  \sum_{n=1}^M (n+1)^b  (1-\frac{1}{n+1})^b  (a_M)^{n+1}-\Sigma_M\\ &
= (M+1)^b  (a_M)^{M+1}- a_M+
\bigO_b\big(\sum_{n=1}^M\frac{1}{n}(n+1)^b  (a_M)^{n+1}\big)\\ &
= (M+1)^b  (a_M)^{M+1}-a_M+\bigO_b((S+1)^{b +1} (a_M)^{S+1})+
\bigO_b\big(\frac{M-S}{S}(M+1)^b  (a_M)^{M+1}\big)\;.
\end{align*}
As $M\ra+\infty$, by taking $S=\frac{4\,M^2}{(1+\sqrt{1+4M}\,)^2}\sim
M$, so that $M-S\sim \sqrt{M}$ and $\frac{M-S}{S}\sim
\frac{1}{\sqrt{M}}$, the sum of the $\bigO_b(\cdot)$ functions in the
above centered line is an $\bigO_b\big(\frac{M^b  (a_M)^M}{\sqrt{M}}\big)$
function. The result follows.
\cqfd

\medskip
In order to simplify the notation in the main body of this text, let 
\begin{equation}\label{defiomegtild}
\wt \omega:t\mapsto
\N_{\E,\omega}(t)= \sum_{x\in\E,\, x\leq t}\omega(x)\;.
\end{equation}
This function is defined on $\RR$ with the usual convention that a sum
over an empty set of indices is $0$. The local finiteness assumption of
the subset $\E$ of $[0,+\infty[$ ensures the finiteness of the growth
function $\wt \omega=\N_{\E,\,\omega}$, and the local finiteness
(hence regularity) of the pair correlation measures $\R^{\F,\psi}_N$
on $\RR$ for $\F= (\,(F_N=\{x\in \E:x\leq N\})_{N\in\NN}, \;\omega)$,
defined in Equation \eqref{eq:defiPCmeasure}. We denote by

\smallskip (PA) the assumption that $\wt \omega(t)\sim c\;t^\alpha\,
e^{\delta t}$ as $t\ra+\infty$ for some constants $c,\delta>0$ and
$\alpha\in\RR$, and by

\smallskip (ET) the assumption that $\wt \omega(t)= c\;t^\alpha\,
e^{\delta t}\big(1+\bigO(e^{-\kappa\,t})\big)$ as $t\ra+\infty$ for
some constants $c,\delta,\kappa>0$ and $\alpha\in\RR$.

\medskip
For all $t\in\RR$ and $\eta>0$, we define the {\it
  $(t,\eta)$-slice of weights} as
\begin{equation}\label{eq:defisliceweight}
\wt\omega(t,\eta)=\wt \omega(t)-\wt \omega(t-\eta)\;.
\end{equation}
The next result describes the asymptotic behaviour for the thin,
though not too thin, slices of weights under one of the two
assumptions (PA) or (ET).

\blemm\label{lem:sliceweight} (1) Let $\eta>0$. Under Assumption (PA),
as $t>0$ tends to $+\infty$, we have
$$
\wt\omega\,(t,\eta)=c\;t^\alpha\,e^{\delta\,t}(1-e^{-\delta\,\eta})
\big(1+\smallo_{\alpha,\eta}(1)\big)\;.
$$

(2) Let $\eta:t\mapsto \eta_t$ be a map from $[0,+\infty[$ to
$]0,1]$. Under Assumption (ET), if ${\displaystyle
\lim_{t\ra+\infty}}\; \eta_t\, e^{\kappa\,t} =+\infty$, then, as $t$
tends to $+\infty$, we have
$$
\wt\omega\,(t,\eta_t)=c\;t^\alpha\,e^{\delta\, t}(1-e^{-\delta\,\eta_t})
\Big(1 +\frac{1}{\eta_t}\,\bigO_{\delta}(e^{-\kappa \,t})+
\bigO_{\alpha,\delta}\big(\frac{1}{t}\big)\Big)\;.
$$
In particular, as $t$ tends to $+\infty$, under these assumptions,
we have
$$
\wt\omega\,(t,\eta_t)=\bigO_{\alpha,\delta}(t^\alpha\,
e^{\delta \,t})\;.
$$
If $\alpha=0$, then as $t$ tends to $+\infty$, we have more precisely
\[
\wt\omega\,(t,\eta_t)=c\;
e^{\delta\, t}(1-e^{-\delta\,\eta_t}) \big(1 +\frac{1}{\eta_t}
\,\bigO_{\delta}(e^{-\kappa \,t})\big)\;.
\]
\elemm

\dem (1) For every $t\in\RR$, let $r_t=c^{-1}t^{-\alpha} \,e^{-\delta\,t}
\;\wt\omega\, (t) -1$, which converges to $0$ as $t\ra +\infty$ since
Assumption (PA) holds. If $t>0$, we have
\begin{align*}
  \wt\omega\,(t,\eta)&=\wt \omega(t)-\wt \omega(t-\eta)
  =c\;t^\alpha\,e^{\delta \,t}(1+r_t)-
  c\;(t-\eta)^\alpha\,e^{\delta \,(t-\eta)}(1+r_{t-\eta})\\ 
  &=c\;t^\alpha\,e^{\delta \,t}\Big(1+r_t
  -\big(1-\frac{\eta}{t}\big)^\alpha e^{-\delta \,\eta}
  (1+r_{t-\eta}))\Big)\\
  &=c\;t^\alpha\,e^{\delta \,t}\Big(1-e^{-\delta \,\eta}+
  r_t-e^{-\delta \,\eta}r_{t-\eta}-e^{-\delta \,\eta}
  \bigO_{\alpha}\big(\frac{\eta}{t}\big)\Big)\;.
\end{align*}
Since $\lim_{t\ra+\infty} \frac{\max\{r_t,\;r_{t-\eta}, \frac{1}{t}\}}
{1-e^{-\delta \,\eta}} =0$, this concludes the proof of Assertion (1).

\medskip
(2) Since $\eta_t\in\;]0,1]$, we have $\bigO(e^{\kappa\,\eta_t})
=\bigO(1)$ as $t\ra+\infty$.  Recall that $(1+s)^\alpha =1+
\bigO_{\alpha}(s)$ as $s\ra 0$. Since Assumption (ET) holds, when
$t>0$ tends to $+\infty$, we hence have
\begin{align*}
  \wt\omega\,(t,\eta_t)&=\wt \omega(t)-\wt \omega(t-\eta_t)
  =c\;t^\alpha\,e^{\delta \,t}\Big(1+\bigO(e^{-\kappa \,t})
  -\big(1-\frac{\eta_t}{t}\big)^\alpha e^{-\delta \,\eta_t}
  (1+\bigO(e^{-\kappa\,t}))\Big)\\
  &=c\;t^\alpha\,e^{\delta \,t}\Big(1-e^{-\delta \,\eta_t}+
  \bigO(e^{-\kappa\,t}) + \bigO_{\alpha}(\frac{\eta_t}{t})\Big)\;.
\end{align*}
Since $\frac{1}{1-e^{-\delta \,\eta_t}}= \bigO_{\delta}(\frac{1}{\eta_t})$ as
$t\ra+\infty$, and since $\lim_{t\ra+\infty}
\frac{e^{-\kappa\,t}}{\eta_t}\, =0$, this proves the result for general
$\alpha$. The proof in the special case when $\alpha=0$ is even
simplier.
\cqfd

\section{An extension of Theorem \ref{theo:mainintro} with error terms}
\label{sec:proof}

We will use in this section the notation $\N_{\E,\,\omega}$ and
$\R^{\F,\psi}_N$ defined in the Introduction, as well as the notation
$\wt \omega(\cdot)$ and $\wt \omega(\cdot,\cdot)$ of Equations
\eqref{defiomegtild} and \eqref{eq:defisliceweight}. For every scaling
function $\psi:\NN\ra[1,+\infty[\,$, we consider the renormalising
    function
\begin{equation}\label{eq:defipsiprim}
\psi':N\mapsto \frac{\N_{\E,\,\omega}(N)^2}{\psi(N)}=
\frac{\wt \omega(N)^{2}}{\psi(N)}\;.
\end{equation}
We denote by $\|\mu\|$ the total mass of a measure $\mu$.

\medskip
We start this section by some comments on the statement of Theorem
\ref{theo:mainintro}.

\medskip\noindent {\bf Remarks. } (1) When $\psi=1$, then
$\psi'(N)=\|\R^{\F,\psi}_N\|$, and the renormalisation in Theorem
\ref{theo:mainintro} (as well as in Theorem \ref{theo:main}) is chosen
in order to obtain probability measures $\frac{1}{\psi'(N)}\,
\R^{\F,\psi}_N$, which turns out to converge to a probability measure
$g_{\F,1}\;dt$ as $N\ra+\infty$.

When $\displaystyle{\lim_{+\infty} \psi} =+\infty$, as the proof below
shows, the renormalisation is precisely chosen in order to obtain a
locally finite nonzero measure, but there is an infinite loss of mass
at infinity, in the sense that ${\displaystyle \lim_{N\ra+\infty}}\;
\frac{1}{\psi'(N)}\, \|\,\R^{\F,\psi}_N\,\|=+\infty$ even though
${\displaystyle \lim_{N\ra+\infty}}\; \frac{1}{\psi'(N)}\,
\|\,(\R^{\F,\psi}_N)_{\mid K}\|$ is finite for every compact subset $K$ of
$[0,+\infty[\,$.

\medskip
(2) The pair correlation measures are sometimes defined (see for
instance \cite{ParPau22a}) by
\[
\wt\R^{\F,\psi}_N=\sum_{x,y\in F_N\;:\;x\neq y}\;\omega(x)\,\omega(y)\,
\Delta_{\psi(N)(y-x)}\;,
\]
that is by adding the assumption $x\neq y$ on the set of pair of
indices $(x,y)\in {F_N}^2$ in the summation (compare with Equation
\eqref{eq:defiPCmeasure}).  Note that when the renormalised measures
$\frac{1}{\psi'(N)}\, \R^{\F,\psi}_N$ weak-star converge to a measure
$\mu$ on $\RR$ which has no atom at $0$ (for instance if the family
$\F$ admits a pair correlation function for the scaling $\psi$ and
renormalisation $\psi'$, as it is the case in Theorem
\ref{theo:mainintro}), then we also have $\lim_{N\ra+\infty}
\frac{1}{\psi'(N)}\; \wt\R^{\F,\psi}_N =\mu$, that is, the
contribution of the diagonal set of indices in the sum defining
$\frac{1}{\psi'(N)}\,\R^{\F,\psi}_N$ is negligible.

\medskip
(3) An archetypical example of a pair $(\E,\omega)$ is given by a
countable set $\wt \E$ endowed with a map $\ell:\wt \E\ra \;
  ]\,0,+\infty\,[\,$ with finite fibers whose image $\E=\ell(\wt \E)$
is locally finite and endowed with the multiplicity function $\omega:
x\mapsto \card (\ell^{-1}(x))$. In this case, we have
\[
\R^{\F,\psi}_N=\sum_{x,y\,\in \,\wt \E\;:\;\ell(x),\, \ell(y)\leq
N} \; \Delta_{\psi(N)(\ell(y)-\ell(x))}\;.
\]
Theorem \ref{theo:mainintro} when $\psi= 1$ then says that if there
exist $\alpha\in\RR$ and $c,\delta>0$ such that, as $t\ra + \infty$,
we have $\card\{x\in\wt \E :\ell(x)\leq t\} \sim c\; t^\alpha \,
e^{\delta \,t}$ then $\frac{1}{\|\R^{\F,1}_N\|} \;\R^{\F,1}_N$
weakstar converges to the measure $\frac{\delta}{2}\;
e^{-\,\delta\,|t|}\,dt$ on the locally compact space $\RR$ as
$N\ra+\infty$.

\bigskip
We now state an extended version with error terms of Theorem
\ref{theo:mainintro}, from which it follows, using Assumption (ET).
Let
\[
g_\delta:t\mapsto\frac{\delta}{2}\; e^{-\,\delta\,|t|}\;.
\]

\btheo\label{theo:main} Let $\E$ be a locally finite subset of
$[0,+\infty[\,$ endowed with a weight function $\omega$, and let
    $\alpha\in\RR$ and $c,\delta>0$. Assume that there exists
    $\kappa>0$ such that as $t\ra+\infty$, we have
\begin{equation}\label{eq:asumpcompET}
\N_{\E,\,\omega}(t)= c\;t^\alpha\,e^{\delta \,t}
(1+\bigO_{\E,\omega}(e^{-\kappa\, t}))\;.
\end{equation}
Let $\psi:\NN\ra[1+\infty[$ be a scaling function and let $A\ge 1$.
For every function $f\in C^1(\RR)$ with compact support in $[-A,A]$,
as $N\ra+\infty$, with Landau functions $\bigO=
\bigO_{\E,\omega,c,\alpha,\delta,\kappa}$ and $\kappa'=\min\{\kappa,
\delta\}$, we have
\[
\frac{\R^{\F,\psi}_N(f)}{\psi'(N)}=\left\{\begin{array}{l}
\int_\RR f(t)\,g_\delta(t)\;dt +
\bigO(A\,e^{\delta A} e^{-\frac{\kappa'}{12} N}\,\|f\|_\infty)
+\bigO(A\,e^{-\frac{\kappa'}{4} N}\,\|f'\|_\infty)
\\{\rm if~} \alpha=0 {\rm ~and~}
\psi= 1, \medskip\\
\int_\RR f(t)\,g_\delta(t)\;dt + \bigO\big(\frac{A^2}{N}\,
(\|f\|_\infty+\|f'\|_\infty)\big)\\{\rm if~} \alpha\neq 0
{\rm ~and~} \psi= 1, \medskip\\
\frac{\delta}{2}\int_\RR f(t)\;dt +
\bigO\big(\frac{A^2\|f\|_\infty}{\psi(N)}\big)+
\bigO\big(A^2\,e^{-\frac{\kappa'}{4} N}\,\psi(N)\,\|f'\|_\infty\big)\\
 {\rm if~}\alpha=0 {\rm ~and~}
\psi {\rm ~converges~to~} +\infty {\rm ~with~subexponential~growth}
\medskip\\
\frac{\delta}{2}\int_\RR f(t)\;dt + 
\bigO\big(\frac{A^2}{\min\{N,\,\psi(N)\}}(\|f\|_\infty+\|f'\|_\infty)\big)
\\{\rm if~}\alpha\neq 0 {\rm ~and~} \psi {\rm ~converges~to~}
+\infty  {\rm ~with~at~most~polynomial~growth.}
\end{array}\right.
\]   
\etheo

\dem Let $\E$, $\omega$, $c$, $\alpha$, $\delta$, $\kappa$ be the
fixed data in the statement of Theorem \ref{theo:main}. Though we
won't indicate the dependency, the Landau functions $\bigO(\cdot)$
below will depend on these fixed data, in the sense defined at the
beginning of Section \ref{sect:prelimposseq}.  Up to replacing
$\kappa$ by $\min\{\kappa,\delta\}$ which does not change the
conclusion of Theorem \ref{theo:main} and is implied by its hypothesis
\eqref{eq:asumpcompET}, we may assume that
\begin{equation}\label{eq:controldeltakappa}
  \kappa\leq\delta\;.
\end{equation}

Let $\psi$, $f$, $A$ and $N$ be the varying data in the statement of
Theorem \ref{theo:main}.  The Landau functions $\bigO(\cdot)$ below
will not depend on these varying data, in the sense defined at the
beginning of Section \ref{sect:prelimposseq}.

Note that if $\iota:t\mapsto -t$, then $\iota_*\R^{\F,\psi}_N=
\R^{\F,\psi}_N$ by using the change of variables $(x,y)\mapsto (y,x)$
in the summation of Equation \eqref{eq:defiPCmeasure}, and that
$g_\delta\circ\iota=g_\delta$. In order to prove Theorem
\ref{theo:main}, we may hence assume by additivity that the support of
$f$ is contained in $[0,A]$, and again by additivity that $f\geq 0$.

Note that $\lim_{+\infty}\psi'=+\infty$ by Equation
\eqref{eq:defipsiprim} since $\wt \omega$ is exponentially growing and
$\psi$ is subexponentially growing in all the cases of Theorem
\ref{theo:main}. As $N\ra+\infty$, by Equations
\eqref{eq:defipsiprim} and \eqref{eq:asumpcompET}, we have
\begin{equation}\label{eq:renormalisation}
  \frac{1}{\psi'(N)}=
  \frac{\psi(N)}{c^2\;N^{2\alpha}\,e^{2\delta N}}
\big(1+\bigO(e^{-\kappa\,N})\big)\,.
\end{equation}
We consider throughout this proof two small quantities $\epsilon,\tau'
\in\;]0,1]$ which will depend on $N$ and
converge to $0$. We define
\[
\tau=\frac{\tau'}{\psi(N)}\;,
\]
and we assume that $\tau\geq 2\epsilon$ when $N$ is large. We will
check this inequality after defining $\epsilon$ and $\tau'$  in
Equations \eqref{eq:defiepsilontauET} and
\eqref{eq:defiepsilontauETalphaneq0}.

\medskip
The following lemma describes the work on the set of indices of some of
the following sums in order to be able to separate the variables $x$
and $y$.

\blemm\label{lem:workonindices}
Let $x,y\in\E$ and $k,n\geq 1$. The system of inequalities
\begin{equation}\label{eq:indexsource}
0<x\leq y\leq N,\;\;\; (k-1)\epsilon< x\leq k\,\epsilon,\;\;\;
(n-1)\tau'< \psi(N)(y-x)\leq n\tau'
\end{equation}
implies the system of inequalities
\begin{equation}\label{eq:indexupperbound}
  (k-1)\epsilon< x\leq k\,\epsilon,\;\;
  (n\tau+k\,\epsilon)-(\tau+\epsilon)< y\leq n\tau+k\,\epsilon,
  \;\;k\leq M^+_\epsilon=
  \Big\lfloor\frac{N-(n-1)\tau}{\epsilon}\Big\rfloor +1\;,
\end{equation}
and is implied by the system of inequalities
\begin{equation}\label{eq:indexlowerbound}
(k-1)\epsilon< x\leq k\,\epsilon,\;\;\;\;
  (n\tau+k\,\epsilon)-\tau< y\leq n\tau+k\,\epsilon-\epsilon,
  \;\;\;\;k\leq M^-_\epsilon=
  \Big\lfloor\frac{N-n\tau}{\epsilon}\Big\rfloor +1\;.
\end{equation}
\elemm

\dem Since we have $\tau=\frac{\tau'}{\psi(N)}$, the last two
inequalities of Equation \eqref{eq:indexsource} are equivalent to
$(n-1)\tau+x< y\leq n\tau +x$. With the middle two inequalities of
Equation \eqref{eq:indexsource}, this implies that
$(n\tau+k\,\epsilon) -(\tau+\epsilon)< y\leq n\tau+k\,\epsilon$ and is
implied by $(n\tau+k\,\epsilon)-\tau< y\leq
n\tau+k\,\epsilon-\epsilon$.

The inequalities $(n\tau+k\,\epsilon) -(\tau+\epsilon)< y$ and $y\leq
N$ imply that $k\,\epsilon\leq N-(n-1)\tau+\epsilon$, hence that
$k\leq M^+_\epsilon$.

The inequalities $y\leq n\tau+k\,\epsilon-\epsilon$ and $k\leq
M^-_\epsilon$ imply that $y\leq n\tau+\epsilon
\big(\frac{N-n\tau}{\epsilon} +1\big)-\epsilon=N$.  The inequalities
$(k-1)\epsilon< x$ and $k\geq 1$ implies that $x>0$. The inequalities
$x\leq k\,\epsilon$, $(n\tau+k\,\epsilon)-\tau< y$ and $n\geq 1$ imply
that $x\leq y$. The result follows.
\cqfd

\medskip
Note that by the definitions of $M_\epsilon^\pm$ in Equations
\eqref{eq:indexupperbound} and \eqref{eq:indexlowerbound}, when
$\tau\geq 2\epsilon$, we have
\begin{equation}\label{eq:epsilonMepsilon}
\epsilon\; M_\epsilon^\pm =
\epsilon\,\Big(\Big\lfloor\frac{N-(n-\frac{1}{2}(1\pm1))\tau}
{\epsilon}\Big\rfloor +1\Big)= N-n\tau+\bigO(\tau)\;.
\end{equation}

\medskip
Let us define
\begin{equation}\label{eq:deinwtmuN}
\mu_N(f)=\sum_{x,y\in \E,\; 0<x\leq y\leq N}
\omega(x)\;\omega(y)\;f(\psi(N)(y-x))\;.
\end{equation}
Since the support of $f$ is contained in $[0,A]$ and $\psi\geq 1$, we
have
\begin{align*}
0\leq \R^{\F,\psi}_N(f)-\mu_N(f)&=\sum_{x,y\in \E,\; 0=x\leq y\leq N}
\omega(x)\;\omega(y)\;f(\psi(N)(y-x))\\ &\leq \sum_{y\in \E\cap[0,A]}
\omega(0)\;\omega(y)\;\|f\|_\infty\;,
\end{align*}
and the term on the right hand is independent of $N$. Note that
$\frac{1}{\psi'(N)}=\bigO(e^{-\kappa N})$ by Equations
\eqref{eq:defipsiprim}, \eqref{eq:asumpcompET} and
\eqref{eq:controldeltakappa}, and by the subexponential growth of
$\psi$ in all the cases of Theorem \ref{theo:main}. Hence in order to
prove Theorem \ref{theo:main}, we therefore only have to prove that
$\frac{\mu_N(f)}{\psi'(N)}$ converges, with the appropriate error
terms, to $\int f\,g_\delta\;dt$ when $\psi= 1$ and to
$\frac{\delta}{2} \int f\;dt$ when $\psi\ra+\infty$.

We will use Riemann sums in order to approximate these integrals, with
subdivision step given by $\tau'$. Since the support of $f$ is
contained in $[0,A]$, we only need to subdivide this interval into
$\big\lceil \frac{A}{\tau'}\big\rceil$ intervals of length $\tau'$.
For all $n\in\{1,\dots,\big\lceil \frac{A}{\tau'}\big\rceil\}$ and
$t\in \;](n-1)\tau',n\tau']$, we have
\begin{equation}\label{eq:acroissementfini}
  f(t)= f(n\tau')-\int_t^{n\tau'}f'(t)\;dt=
  f(n\tau')+\bigO\big(\int_{(n-1)\tau'}^{n\tau'}|f'(t)|\;dt\big)
  =f(n\tau')+\bigO(\tau'\,\|f'\|_\infty)\;.
\end{equation}
Since $f\geq 0$, we may assume that
$f(n\tau')+\bigO(\tau'\,\|f'\|_\infty)\geq 0$. Let us define
\begin{equation}\label{eq:definanN}
  a_{n,N}=\sum_{k=1}^{+\infty}\sum_{\substack{x,y\in\E\;:\;0<x\leq y\leq N,\\
    (k-1)\epsilon< x\leq k\,\epsilon,\\
(n-1)\tau'< \psi(N)(y-x)\leq n\tau'}} \omega(x)\;\omega(y)
\end{equation}
which is a nonnegative finite sum and depends also on $\epsilon$ and
on $\tau'$. With the help of Equation \eqref{eq:acroissementfini},
Equation \eqref{eq:deinwtmuN} becomes, by subdividing the range of $x$
into half-open intervals of length $\epsilon$ and the range of
$\psi(N)(y-x)$ (with values at most $A$ in order to be in the support
of $f$) into half-open intervals of length $\tau'$,
\begin{align}
  \mu_N(f)&
  =\sum_{n=1}^{\lceil \frac{A}{\tau'}\rceil}\sum_{k=1}^{+\infty}
  \sum_{\substack{x,y\in\E\;:\;0<x\leq y\leq N,\\
      (k-1)\epsilon< x\leq k\,\epsilon,\\
      (n-1)\tau'< \psi(N)(y-x)\leq n\tau'}}
  \omega(x)\;\omega(y)\;f(\psi(N)(y-x))
  \nonumber\\&=\sum_{n=1}^{\lceil \frac{A}{\tau'}\rceil}
  \big(f(n\tau')+\bigO(\tau'\,\|f'\|_\infty)\big)\;a_{n,N}\;.
  \label{eq:relatmuNfwithanN}
\end{align}
We now proceed by a majoration and minoration of $a_{n,N}$. Let us define
\begin{equation}\label{eq:definanNplus}
  a_{n,N}^+=\sum_{k=1}^{M_\epsilon^+}\sum_{\substack{x,y\in\E\;:\;
      (k-1)\epsilon< x\leq k\,\epsilon\\
      (n\tau+k\,\epsilon)-(\tau+\epsilon)< y\leq n\tau+k\,\epsilon}}
  \omega(x)\;\omega(y)
\end{equation}
and
\begin{equation}\label{eq:definanNminus}
  a_{n,N}^-=\sum_{k=1}^{M_\epsilon^-}\sum_{\substack{x,y\in\E\;:\;
      (k-1)\epsilon< x\leq k\,\epsilon\\
      (n\tau+k\,\epsilon)-\tau< y\leq (n\tau+k\,\epsilon)+\epsilon}}
  \omega(x)\;\omega(y)\;.
\end{equation}
By Lemma \ref{lem:workonindices}, we have
\[
a_{n,N}^-\leq a_{n,N}\leq a_{n,N}^+\;.
\]
Since the variables $x,y$ are separated in the sums defining
$a_{n,N}^\pm$ and by the definition \eqref{eq:defisliceweight} of the
slices of weights, we have
\begin{equation}\label{eq:definanNplussliced}
  a_{n,N}^+= \sum_{k=1}^{M_\epsilon^+}\;\wt\omega\,(k\epsilon,\epsilon)\;
  \wt\omega\,(n\tau+k\epsilon,\tau+\epsilon)
\end{equation}
and
\begin{equation}\label{eq:defianNminussliced}
  a_{n,N}^-=\sum_{k=1}^{M_\epsilon^-}\;\wt\omega\,(k\epsilon,\epsilon)\;
  \wt\omega\,(n\tau+k\epsilon+\epsilon,\tau-\epsilon)\;.
\end{equation}
We study the quantities $a_{n,N}^\pm$ under the assumption
\eqref{eq:asumpcompET} on the asymptotic behaviour of the weights.

\medskip
\noindent{\bf Asymptotics on $a_{n,N}^\pm$.} By Equation
\eqref{eq:definanNplussliced} and by two applications of Lemma
\ref{lem:sliceweight} (2) with $(t,\eta_t)=(k\epsilon,\epsilon)$ and
$(t,\eta_t)=(n\tau+ k\epsilon,\tau +\epsilon)$ as $k\ra+\infty$, and
up to verifying when we will define $\epsilon$ and $\tau'$ that the
assumption of this lemma (besides Assumption (ET) which holds by
Equation \eqref{eq:asumpcompET}) is satisfied, for all $N$ large
enough and $n\in\{1,\dots, \lceil\frac{A}{\tau'}\rceil\}$, we have
\begin{align}
  a_{n,N}^+= \sum_{k=1}^{M_\epsilon^+}\; &
  c\;(k\epsilon)^\alpha\,e^{\delta k\epsilon}(1-e^{-\delta\epsilon})
\Big(1 +\frac{1}{\epsilon}\,\bigO(e^{-\kappa k\epsilon})+
\bigO\big(\frac{1}{k\epsilon}\big)\Big)\nonumber\\ 
\times \;&c\;(n\tau+ k\epsilon)^\alpha\,e^{\delta(n\tau+ k\epsilon)}
(1-e^{-\delta(\tau+ \epsilon)})\Big(1 +\frac{1}{\tau+ \epsilon}\,
\bigO(e^{-\kappa(n\tau+ k\epsilon})+
\bigO\big(\frac{1}{n\tau+ k\epsilon}\big)\Big)\nonumber\\ = c^2 &
(1-e^{-\delta\epsilon})(1-e^{-\delta(\tau+ \epsilon)})\;e^{\delta n\tau}
\sum_{k=1}^{M_\epsilon^+}\; z^+_k\;,\label{eq:relatanNzkplus}
\end{align}
where{\small
\[
z^+_k=(k\epsilon)^{2\alpha}\,e^{2\delta k\epsilon}\,
\big(1+ \frac{n\tau}{k\epsilon}\big)^\alpha
\Big(1 +\frac{1}{\epsilon}\,\bigO(e^{-\kappa k\epsilon})+
\bigO\big(\frac{1}{k\epsilon}\big)\Big)
\Big(1 +\frac{1}{\tau+ \epsilon}\,\bigO(e^{-\kappa(n\tau+ k\epsilon)})+
\bigO\big(\frac{1}{n\tau+ k\epsilon}\big)\Big)\;.
\]}
$\!\!$Since $\tau\pm\epsilon\geq \epsilon$, $n\tau+k\epsilon\geq k\epsilon$
and $e^{-\kappa n\tau} \leq 1$, this simplifies as
\begin{equation}\label{eq:defizkplus}
z^+_k=(k\epsilon)^{2\alpha}\,e^{2\delta k\epsilon}\,
\big(1+ \frac{n\tau}{k\epsilon}\big)^\alpha
\Big(1 +\frac{1}{\epsilon}\,\bigO(e^{-\kappa k\epsilon})+
\bigO\big(\frac{1}{k\epsilon}\big)\Big)^2\;.
\end{equation}

\medskip
\noindent{\bf Case 1: Let us first assume that $\alpha=0$.}  We define
in this case
\begin{equation}\label{eq:defiepsilontauET}
  \epsilon=e^{-\frac{\kappa}{3}N} {\rm ~~~and~~~}
  \tau'= e^{-\frac{\kappa}{4}N}\psi(N)\;.
\end{equation}
In particular, we have $\lim_{N\ra+\infty} \epsilon\,e^{\kappa N}
=+\infty$, we have $\tau\geq 2\epsilon$ if $N$ is large enough, and
since $\psi$ grows subexponentially under the two assumptions of
Theorem \ref{theo:main} when $\alpha=0$, the quantity $\tau'$ tends to
$0$ as $N\ra+\infty$.  By the last claim of Lemma
\ref{lem:sliceweight} (2), Equation \eqref{eq:defizkplus} for $z^+_k$
simplifies as
\begin{align*}
z^+_k=e^{2\delta k\epsilon}
\big(1 +\frac{1}{\epsilon}\,\bigO(e^{-\kappa k\epsilon})\big)^2
= e^{2\delta k\epsilon}+\frac{1}{\epsilon^2}
\,\bigO(e^{(2\delta-\kappa) k\epsilon})\;.
\end{align*}
Hence by a geometric series summation, since $2\delta> \delta\geq
\kappa$ by Equation \eqref{eq:controldeltakappa}, we have
\[
\sum_{k=1}^{M_\epsilon^+}\; z^+_k=
\frac{e^{2\delta \epsilon(M_\epsilon^++1)}-1}{e^{2\delta \epsilon}-1}+
\frac{1}{\epsilon^2}\bigO\big(
\,\frac{e^{(2\delta-\kappa)\epsilon(M_\epsilon^++1)}-1}
       {e^{(2\delta-\kappa) \epsilon}-1}\big)\;.
\]
Note that $\frac{e^{2\delta \epsilon}-1} {e^{(2\delta-\kappa)
    \epsilon}-1}=\bigO(1)$ as $\epsilon\ra 0$, and recall that
$n\tau=\frac{n\tau'}{\psi(N)}\leq \frac{A+1}{\psi(N)}$ for every
$n\in\{1,\dots, \lceil\frac{A}{\tau'}\rceil\}$. By Equation
\eqref{eq:epsilonMepsilon}, by the definition
\eqref{eq:defiepsilontauET} of $\epsilon$ and since
$2\delta\geq\delta\geq \kappa\geq \frac{\kappa}{3}$ by Equation
\eqref{eq:controldeltakappa}, we hence have
\begin{align*}
\sum_{k=1}^{M_\epsilon^+}\; z^+_k&=
\frac{e^{2\delta(N -n\tau)}}{e^{2\delta \epsilon}-1}\big(e^{\bigO(\tau)}-
e^{\delta(n\tau-N)}+\frac{1}{\epsilon^2}\bigO(e^{\kappa(n\tau-N)})\big)
\\ &=\frac{e^{2\delta(N -n\tau)}}{e^{2\delta \epsilon}-1}
\big((1+\bigO(e^{-\frac{\kappa}{4}N}))+\bigO(e^{\frac{\delta A}{\psi(N)}}
e^{-\delta N})+\bigO(e^{\frac{\kappa A}{\psi(N)}}
e^{-\frac{\kappa}{3}N})\big)\\ &=
\frac{e^{2\delta(N -n\tau)}}{e^{2\delta \epsilon}-1}
\big(1+\bigO(e^{\frac{\delta A}{\psi(N)}} e^{-\frac{\kappa}{4} N})\big)\;.
\end{align*}
By the definition \eqref{eq:defiepsilontauET} of $\epsilon$ and
$\tau'$ (so that $\tau=e^{-\frac{\kappa}{4}N}$), we have
\[
1-e^{-\delta(\tau\pm\epsilon)}=
\delta\tau(1\pm\frac{\epsilon}{\tau})(1+\bigO(\tau+\epsilon))=
\delta\tau(1+\bigO(e^{-\frac{\kappa}{12}N}))\;,
\]
and
\[
\frac{1-e^{-\delta\epsilon}}{e^{2\delta\epsilon}-1}=
\frac{1}{2}(1+\bigO(\epsilon))=
\frac{1}{2}(1+\bigO(e^{-\frac{\kappa}{3} N}))\;.
\]
Therefore Equation \eqref{eq:relatanNzkplus} becomes
\begin{align*}
a_{n,N}^+&=c^2\;e^{2\delta N}
\frac{1-e^{-\delta\epsilon}}{e^{2\delta\epsilon}-1}
(1-e^{-\delta(\tau+\epsilon)})\,e^{-\delta n\tau}\,
\big(1+\bigO(e^{\frac{\delta A}{\psi(N)}} e^{-\frac{\kappa}{4} N})\big)\\&
=(c^2\;e^{2\delta N})
(\,\frac{\delta}{2}\;\tau\;e^{-\delta n\tau})
\big(1+\bigO(e^{\frac{\delta A}{\psi(N)}} e^{-\frac{\kappa}{12} N})\big)\;.
\end{align*}
Taking into account the small differences between $a_{n,N}^+$ and
$a_{n,N}^-$ in Equations \eqref{eq:definanNplus} and
\eqref{eq:definanNminus}, a similar computation gives the same formula
for $a_{n,N}^-$. Since $a_{n,N}^-\leq a_{n,N}\leq a_{n,N}^+$, we hence
have
\begin{equation}\label{eq:asymptotanN}
  a_{n,N}=(c^2\;e^{2\delta N})
(\,\frac{\delta}{2}\;\tau\;e^{-\delta n\tau})
\big(1+\bigO(e^{\frac{\delta A}{\psi(N)}} e^{-\frac{\kappa}{12} N})\big)\;.
\end{equation}

\medskip
\noindent{\bf End of the proof of Theorem \ref{theo:main} when
  $\alpha=0$.}  Note that $\sum_{n=1}^{\lceil\frac{A}{\tau'}\rceil}
e^{-\frac{\delta n\tau'}{\psi(N)}}=\bigO(\frac{A}{\tau'})$ as
$N\ra+\infty$. Since $\tau= \frac{\tau'}{\psi(N)}$, by Equations
\eqref{eq:relatmuNfwithanN}, \eqref{eq:renormalisation} and
\eqref{eq:asymptotanN}, for $N$ large enough, we have
\begin{align}
 \frac{\mu_N(f)}{\psi'(N)}&= \sum_{n=1}^{\lceil\frac{A}{\tau'}\rceil}
 \big(f(n\tau')+\bigO(\tau'\,\|f'\|_\infty)\big)
 \,\frac{\delta}{2}\,\tau'\,e^{-\frac{\delta n\tau'}{\psi(N)}}
\big(1+\bigO(e^{\frac{\delta A}{\psi(N)}} e^{-\frac{\kappa}{12} N})\big)
\big(1+\bigO(e^{-\kappa\,N})\big)
\nonumber\\ & =\Big(\sum_{n=1}^{\lceil\frac{A}{\tau'}\rceil}
f(n\tau')\,\frac{\delta}{2}\,\tau'\,e^{-\frac{\delta n\tau'}{\psi(N)}}\Big)
+\bigO(A\,e^{\frac{\delta A}{\psi(N)}} e^{-\frac{\kappa}{12} N}\|f\|_\infty)
+\bigO(A\,\tau'\,\|f'\|_\infty)\;.\label{eq:globalestimate1}
\end{align}

Assume first that $\psi= 1$. Recall that $g_\delta:t\mapsto
\frac{\delta}{2}\,e^{-\delta t}$ is bounded with bounded derivative on
$[0,+\infty[\,$. By the standard Riemann sum approximation with error
term of an integral, and since the support of $f$ is contained in
$[0,A]$ with $A\geq 1$, we have
\begin{align}
\sum_{n=1}^{\lceil\frac{A}{\tau'}\rceil}
f(n\tau')\,\frac{\delta}{2}\,\tau'\,e^{-\delta n\tau'}&=
\int_0^{+\infty}f(t)\,g_\delta(t)\;dt +
\bigO\big(\tau'(\|fg_\delta\|_\infty+\Var(fg_\delta))\big)\nonumber\\& =
\int_0^{+\infty}f(t)\,g_\delta(t)\;dt +
\bigO\big(A\,\tau'(\|f\|_\infty+\|f'\|_\infty)\big)\;.
\label{eq:Riemansumestim1}
\end{align}
With Equation \eqref{eq:globalestimate1}, this proves Theorem
\ref{theo:main} when $\alpha=0$ and $\psi= 1$.

\medskip
Assume now that $\displaystyle{\lim_{+\infty}\psi}=+\infty$. Note that
$e^{-\frac{\delta n\tau'}{\psi(N)}}=1+\bigO(\frac{A}{\psi(N)})$ since
$n\leq \lceil\frac{A}{\tau'}\rceil$.  Hence a similar Riemann sum
argument gives
\begin{align}
\sum_{n=1}^{\lceil\frac{A}{\tau'}\rceil}
f(n\tau')\,\frac{\delta}{2}\,\tau'\,e^{-\delta n\frac{\tau'}{\psi(N)}}&=
\Big(\,\frac{\delta}{2}\int_0^{+\infty}f(t)\;dt +
\bigO\big(A\,\tau'(\|f\|_\infty+\|f'\|_\infty)\big)\Big)
(1+\bigO(\frac{A}{\psi(N)}))\nonumber\\ &=
\frac{\delta}{2}\int_0^{+\infty}f(t)\;dt+
\bigO\big(\frac{A^2\|f\|_\infty}{\psi(N)}\big)+
\bigO\big(A^2\,\tau'(\|f\|_\infty+\|f'\|_\infty)\big)\;.
\label{eq:Riemansumestim2}
\end{align}
Since $\psi$ grows subexponentially, Equations
\eqref{eq:globalestimate1} and \eqref{eq:Riemansumestim2} imply
Theorem \ref{theo:main} when $\alpha=0$ and
$\displaystyle{\lim_{+\infty}\psi}=+\infty$.

\medskip
\noindent{\bf Case 2: Let us now assume that $\alpha\neq 0$.} The
scheme of proof is the same one as in Case 1, though more
technical. We will only be able to obtain the result under a bit
stronger assumption on the scaling function $\psi$ and with a much
weaker error term, due to the polonomial term in the asymptotic growth
of the counting function $\wt\omega$. Since $\psi$ is assumed to have
polynomial growth, we fix $\ga'\geq 1$ such that $\lim_{N\ra+\infty}
\frac{\psi(N)}{N^{\ga'-1}}=0$. We define in this case
\begin{equation}\label{eq:defiepsilontauETalphaneq0}
  \epsilon=\frac{1}{N^{2\ga'}} {\rm ~~~and~~~}
  \tau'= \frac{1}{N^{\ga'}}\,\psi(N)\;.
\end{equation}
In particular, we have ${\displaystyle \lim_{N\ra+\infty}}\; \epsilon \,e^{\kappa N} =
+\infty$, we have $\tau\geq 2\epsilon$ if $N$ is large enough, and
$\tau'$ tends to $0$ as $N\ra+\infty$.

For $N$ is large enough, since $n\tau=\frac{n\tau'}{\psi(N)}$ is
bounded by $A+1\leq 2A$ for $n\leq \lceil\frac{A}{\tau'}\rceil$ and
since $A\geq 1$, as $k\ra+\infty$, Equation \eqref{eq:defizkplus}
gives
\begin{align*}
z^+_k&=(k\epsilon)^{2\alpha}\,e^{2\delta k\epsilon}\,
\big(1+ \bigO(\frac{n\tau}{k\epsilon})\big)
\Big(1 +\frac{1}{\epsilon}\,\bigO(e^{-\kappa k\epsilon})+
\bigO\big(\frac{1}{k\epsilon}\big)\Big)^2\nonumber\\&=
(k\epsilon)^{2\alpha}\,e^{2\delta k\epsilon}+
\bigO\big(A\,(k\epsilon)^{2\alpha-1}\,e^{2\delta k\epsilon}\big)+
\frac{1}{\epsilon^2}\,\bigO\big(A\,(k\epsilon)^{2\alpha}\,
e^{(2\delta-\kappa) k\epsilon}\big)\;.\label{eq:defizkplusalphneq0}
\end{align*}
We now apply Lemma \ref{lem:gengeomlimit} with $M=M_\epsilon^+$, with
$b= 2\alpha$ or $b= 2\alpha-1$, and with $a_M= e^{2\delta \epsilon}$
or $a_M= e^{(2\delta-\kappa) \epsilon}$. The hypothesis of this lemma
is satisfied, since by the definition of $\epsilon$ in Equation
\eqref{eq:defiepsilontauETalphaneq0} and of $M_\epsilon^\pm$ in
Equations \eqref{eq:indexupperbound} and \eqref{eq:indexlowerbound},
and with $b'=2\delta$ or $b'=2\delta-\kappa$ for every
$\ga\in\;]\frac{2\ga'}{2\ga'+1},1[\,$, we have
\[
(M_\epsilon^\pm)^\ga \ln (e^{b' \epsilon})\sim b'\,N^\ga \,\epsilon^{1-\ga}=
b'\,N^{\ga-(1-\ga)(2\ga')}=b'\,N^{\ga(2\ga'+1)-2\ga'}\;,
\]
which converges to $+\infty$ as $N\ra+\infty$ by the assumption on
$\ga$. We hence have
\begin{align*}
  \sum_{k=1}^{M_\epsilon^+}\; z^+_k=\;&
  \frac{e^{2\delta\epsilon}} {e^{2\delta \epsilon}-1}\,
  (\epsilon M_\epsilon^+)^{2\alpha} \,e^{2\delta \,\epsilon M_\epsilon^+}
  \big(1+\bigO(\frac{\sqrt{\epsilon}}{\sqrt{\epsilon M_\epsilon^+}})\big)
  +\bigO\big(A\frac{e^{2\delta \epsilon}}{e^{2\delta \epsilon}-1}\,
  (\epsilon M_\epsilon^+)^{2\alpha-1} e^{2\delta \epsilon M_\epsilon^+}\big)
  \\&
  +\frac{1}{\epsilon^2}\bigO\big(A\frac{e^{(2\delta-\kappa) \epsilon}}
    {e^{(2\delta-\kappa)\epsilon}-1}\,(\epsilon M_\epsilon^+)^{2\alpha}\,
    e^{(2\delta-\kappa) \epsilon M_\epsilon^+}\big)\;.
\end{align*}
Since $\frac{e^{2\delta \epsilon}-1}{e^{(2\delta-\kappa)\epsilon}-1}=
\bigO(1)$ as $\epsilon\ra 0$, by Equation \eqref{eq:epsilonMepsilon},
since $n\tau\leq 2A$ as seen above, and since $\sqrt{\epsilon}=
\frac{1}{N^{\ga'}}\leq \frac{1}{\sqrt{N}}$, for $N$ large enough, we
have
\begin{align*}
  \sum_{k=1}^{M_\epsilon^+}\; z^+_k=\;&
\frac{N^{2\alpha} e^{2\delta N}}{e^{2\delta \epsilon}-1}\,e^{-2\delta n\tau}
\Big(e^{\bigO(\epsilon)}(1-\frac{n\tau+\bigO(\tau)}{N})^{2\alpha}
e^{\bigO(\tau)} \big(1+\bigO(\frac{\sqrt{\epsilon}}{\sqrt{N}})\big)
\\ & +\bigO(\frac{A}{N})+ \frac{1}{\epsilon^2}\bigO(A\,e^{-\kappa N})
\Big) 
=\frac{N^{2\alpha} e^{2\delta N}}{e^{2\delta \epsilon}-1}\,
e^{-2\delta n\tau}\big(1+\bigO(\frac{A}{N})\big)\;.
\end{align*}
By the definition \eqref{eq:defiepsilontauETalphaneq0} of $\epsilon$
and $\tau'$ (so that $\tau=\frac{1}{N^{\ga'}}$), we have
\[
1-e^{-\delta(\tau\pm\epsilon)}=
\delta\tau(1\pm\frac{\epsilon}{\tau})(1+\bigO(\tau+\epsilon))=
\delta\tau\big(1+\bigO(\frac{1}{N^{\ga'}})\big)\;,
\]
and
\[
\frac{1-e^{-\delta\epsilon}}{e^{2\delta\epsilon}-1}=
\frac{1}{2}(1+\bigO(\epsilon))=
\frac{1}{2}\big(1+\bigO(\frac{1}{N^{2\ga'}})\big)\;.
\]
Since $\ga'\geq 1$, Equation \eqref{eq:relatanNzkplus} therefore
becomes
\begin{align*}
a_{n,N}^+&=c^2\;N^{2\alpha}\;e^{2\delta N}
\frac{1-e^{-\delta\epsilon}}{e^{2\delta\epsilon}-1}
(1-e^{-\delta(\tau+\epsilon)})\,e^{-\delta n\tau}\,
\big(1+\bigO(\frac{A}{N})\big)\\&
=(c^2\;N^{2\alpha}\;e^{2\delta N})
(\,\frac{\delta}{2}\;\tau\;e^{-\delta n\tau})
\big(1+\bigO(\frac{A}{N})\big)\;.
\end{align*}
As when $\alpha=0$, computing similarly $a_{n,N}^-$, we have
\begin{equation}\label{eq:asymptotanNalphaneq0}
  a_{n,N}=(c^2\;N^{2\alpha}\;e^{2\delta N})
  (\,\frac{\delta}{2}\;\tau\;e^{-\delta n\tau})
\big(1+\bigO(\frac{A}{N})\big)\;.
\end{equation}

\medskip
\noindent{\bf End of the proof of Theorem \ref{theo:main} when
  $\alpha\neq 0$. } As when $\alpha=0$, by Equations
\eqref{eq:relatmuNfwithanN}, \eqref{eq:renormalisation} and
\eqref{eq:asymptotanNalphaneq0}, for $N$ large enough, we have
\begin{align}
 \frac{\mu_N(f)}{\psi'(N)}&= \sum_{n=1}^{\lceil\frac{A}{\tau'}\rceil}
 \big(f(n\tau')+\bigO(\tau'\,\|f'\|_\infty)\big)
 \,\frac{\delta}{2}\,\tau'\,e^{-\frac{\delta n\tau'}{\psi(N)}}
\big(1+\bigO(\frac{A}{N})\big)
\big(1+\bigO(e^{-\kappa\,N})\big)
\nonumber\\ & =\Big(\sum_{n=1}^{\lceil\frac{A}{\tau'}\rceil}
f(n\tau')\,\frac{\delta}{2}\,\tau'\,e^{-\frac{\delta n\tau'}{\psi(N)}}\Big)
+\bigO(\frac{A^2}{N}\|f\|_\infty)
+\bigO(A\,\tau'\,\|f'\|_\infty)\;.\label{eq:globalestimate2}
\end{align}

When $\psi=1$, since $\tau'=N^{-\ga'}$ with $\ga'\geq 1$ and
$A\geq 1$, Equations \eqref{eq:globalestimate2} and
\eqref{eq:Riemansumestim1} prove Theorem \ref{theo:main} when
$\alpha\neq0$ and $\psi= 1$.

When $\displaystyle{\lim_{+\infty}\psi}=+\infty$, since
$\lim_{N\ra+\infty}\frac{\psi(N)}{N^{\ga'-1}}=0$ so that
\[
\tau'=\frac{\psi(N)}{N^{\ga'}}=\frac{1}{N}\frac{\psi(N)}{N^{\ga'-1}}=
\bigO\big(\frac{1}{\min\{N,\psi(N)\}}\big)\;,
\]
Equations \eqref{eq:globalestimate2} and \eqref{eq:Riemansumestim2}
imply Theorem \ref{theo:main} when $\alpha\neq 0$ and
$\displaystyle{\lim_{+\infty}\psi}=+\infty$.  \cqfd

\bigskip
Let us now prove the analogous result under Assumption (PA) when
$\psi=1$, following closely the scheme of proof of Theorem
\ref{theo:main}, and using the same notation. This is useful in order
to deal with counting asymptotics that sometimes do not come with a
known error term.

\btheo\label{theo:mainPA} Let $\E$ be a locally finite subset of
$[0,+\infty[\,$ endowed with a weight function $\omega$, and let
$\alpha\in\RR$ and $c,\delta>0$. Assume that as $t\ra+\infty$, we
have
\[
\N_{\E,\,\omega}(t)\sim c\;t^\alpha\,e^{\delta \,t}\;.
\]
Then the family $\F=(\,(F_N=\{x\in \E:x\leq N\})_{N\in\NN}, \;\omega)$
admits a pair correlation function for the scaling function $\psi=1$
and renormalizing function $\psi':N\mapsto \N_{\E,\,\omega}(N)^2$,
which is equal to $g_\delta:t\mapsto \frac{\delta}{2}\,e^{-\delta|t|}$.
\etheo

\dem 
In order to prove, as requested for the weak-star convergence, that
$\frac{1}{\psi'(N)}\,\R^{\F,\psi}_N(f)$ converges to $\int_\RR
f\,g_\delta\;dt$ for every continuous fonction $f$ with compact support
on $\RR$, we may assume by density that $f$ is of class $C^1$.

Since $\psi=1$, analogously with the beginning of the proof of Theorem
\ref{theo:main}, as $N\ra+\infty$, we have
\begin{equation}\label{eq:renormalisationPA}
 \frac{1}{\psi'(N)}=
  \frac{1}{c^2\;N^{2\alpha}\,e^{2\delta N}}(1+\smallo(1))\,.
\end{equation}
The first part of the proof is identical with the proof of Theorem
\ref{theo:main} until Equation \eqref{eq:defianNminussliced}, and we
will not repeat it here. In the same way, we are lead to study the
quantities $a_{n,N}^\pm$ given by the equations
\eqref{eq:definanNplussliced} and \eqref{eq:defianNminussliced}. Here,
we do not separate the treatment of the cases $\alpha=0$ and
$\alpha\ne 0$.

\medskip
\noindent{\bf Asymptotics on $a_{n,N}^\pm$.}  Let
$\sigma_\pm=\frac{1-\pm1}{2}$. Note that
$(k\epsilon)^{2\alpha}e^{2\delta k\epsilon}$ tends to $+\infty$ as
$k\ra+\infty$, in order to control the beginning of the following
summations over $k$. Furthermore, since $\epsilon\leq 1\leq A$ and
$n\tau\leq n\tau'\leq A+1$, as $k\ra+\infty$, we have $(1+\frac{n\tau
  +\sigma_\pm\epsilon} {k\epsilon} )^\alpha = 1+
\bigO(\frac{A}{k\epsilon})$.  In the following estimates, except for
$\bigO(\tau)$ taken as $\tau\ra0$ uniformly on everything else unless
indicated by a subscript, the Landau functions $\bigO$ and $\smallo$
are taken as $N\ra+\infty$, uniformly in $n\in\{1,\dots,\lceil
\frac{A}{\tau'}\rceil\}$ up to Equation \eqref{eq:controlanNpm}. By
Equations \eqref{eq:definanNplussliced} and
\eqref{eq:defianNminussliced}, we hence have
\begin{align}
a_{n,N}^\pm=\;& \Big(\sum_{k=1}^{M_\epsilon^\pm}\;c\;(k\epsilon)^\alpha\,
e^{\delta k\epsilon}(1-e^{-\delta\epsilon})\nonumber\\ &\;\;\;\;\;\;\;\; \times
c\;(n\tau+ k\epsilon+\sigma_\pm\epsilon))^\alpha \,e^{\delta(n\tau+
  k\epsilon+\sigma_\pm\epsilon)}(1-e^{-\delta(\tau\pm \epsilon)})\Big)
(1+\smallo_{\epsilon,\tau}(1))\nonumber\\ =\;&
c^2(1-e^{-\delta\epsilon})(1-e^{-\delta(\tau\pm \epsilon)})\;e^{\delta n\tau}
\Big(\sum_{k=1}^{M_\epsilon^+}(k\epsilon)^{2\alpha}e^{2\delta k\epsilon}\Big)
(1+\smallo_{\epsilon,\tau,A}(1))\,\;\;.\label{eq:anNhypPA}
\end{align}
By Lemma \ref{lem:gengeomlimit} with $M=M^\pm_\epsilon$ which goes to
$+\infty$ as $N\ra+\infty$ when $\epsilon$ is fixed, and
$a_M=e^{2\delta \epsilon}$ which is constant when $\epsilon$ is fixed, by
Equation \eqref{eq:epsilonMepsilon}, since $\epsilon\leq
\frac{\tau}{2}$ and since $(1-\frac{n\tau +\bigO(\tau)}{N})^\alpha =
1+\bigO(\frac{A}{N})$, as $N\ra+\infty$, we have
\begin{align*}
  \sum_{k=1}^{M_\epsilon^+}(k\epsilon)^{2\alpha}e^{2\delta k\epsilon}&=
  \frac{e^{2\delta \epsilon}}{e^{2\delta \epsilon}-1}\,
  (\epsilon M_\epsilon^+)^{2\alpha} \,e^{2\delta \,\epsilon M_\epsilon^+}
  (1+\smallo_{\epsilon,\tau}(1))\\ & =
  \frac{e^{\bigO(\tau)}}{e^{2\delta \epsilon}-1}\,
  N^{2\alpha}\,e^{2\delta (N-n\tau)}(1+\smallo_{\epsilon,\tau,A}(1))\;.
\end{align*}
Therefore Equation \eqref{eq:anNhypPA} becomes
\begin{equation}\label{eq:controlanNpm}
a_{n,N}^\pm=c^2\,N^{2\alpha}\,e^{2\delta N}\,e^{\bigO(\tau)}\,
\frac{1-e^{-\delta\epsilon}}{e^{2\delta \epsilon}-1}
(1-e^{-\delta(\tau\pm \epsilon)})\;
e^{-\delta n\tau}(1+\smallo_{\epsilon,\tau,A}(1))\;.
\end{equation}

\medskip
\noindent{\bf End of the proof. } By Equations
\eqref{eq:relatmuNfwithanN} and \eqref{eq:renormalisationPA}, since
$a_{n,N}^-\leq a_{n,N}\leq a_{n,N}^+$ and $\tau'=\tau$, we have
\[
\limsup_{N\ra+\infty} \frac{\mu_N(f)}{\psi'(N)}\leq
\sum_{n=1}^{\lceil\frac{A}{\tau'}\rceil}
\frac{1-e^{-\delta\epsilon}}{e^{2\delta\epsilon}-1}
e^{\bigO(\tau')}\,\big(f(n\tau')+\bigO_f(\tau')\big)\;
(1-e^{-\delta(\tau'+\epsilon)})\,e^{-\delta n\tau'}\;.
\]
By taking the limit as $\epsilon\ra0$, we then have
\[
\limsup_{N\ra+\infty} \frac{\mu_N(f)}{\psi'(N)}\leq
\sum_{n=1}^{\lceil\frac{A}{\tau'}\rceil}
e^{\bigO(\tau')}\,\big(f(n\tau')+\bigO_f(\tau')\big)\;
\frac{1-e^{-\delta\tau'}}{2}\,e^{-\delta n\tau'}\;.
\]
Since $1-e^{-\delta\tau'}\sim \delta\tau'$ as $\tau'\ra0$, by taking
the limit as $\tau'\ra0$ and by a Riemann sum argument, we have
\[
\limsup_{N\ra+\infty} \frac{\mu_N(f)}{\psi'(N)}\leq
\int_0^A f(t)\,g_\delta(t)\;dt\;.
\]
A similar computation gives
\[
\liminf_{N\ra+\infty} \frac{\mu_N(f)}{\psi'(N)}\geq
\int_0^A f(t)\,g_\delta(t)\;dt\;,
\]
which proves Theorem \ref{theo:mainPA}.
\cqfd

\section{Geometric applications}
\label{sec:application}

In this section, we apply the Theorems \ref{theo:mainPA} and
\ref{theo:main} to the sets (with multiplicities) of the lengths of
closed geodesics and common perpendiculars in negatively curved
spaces, and to other discrete sets with similar growth properties that
arise in geometry and dynamics.  We assume familiarity with geometry
and ergodic theory in negative curvature, and we refer, for instance,
to \cite{BroParPau19} for more background and for definitions of the
various objects below.

Let $X$ be either a proper $\RR$-tree without terminal points or a
complete simply connected Riemannian manifold with pinched negative
curvature at most $-1$. Let $\Ga$ be a nonelementary discrete group of
isometries of $X$. Assume that the critical exponent $\delta_\Ga$ of
$\Ga$ is finite.  Let $\D^\pm=(D^\pm_k)_{k\in I^\pm}$ be locally
finite $\Ga$-equivariant families of nonempty proper closed convex
subsets of $X$.  Assume that the outer and inner skinning measures
$\sigma^\pm_{\D^\mp}$ of the families $\D^\mp$ are finite and nonzero,
and that the Bowen-Margulis measure $m_{\rm BM}$ is finite and mixing
for the geodesic flow on the space $\Ga\bs\G X$ of geodesic lines of
$X$ modulo $\Ga$.

A {\em common perpendicular} from $\pi(D^-_k)$ to $\pi(D^+_j)$ is a
locally geodesic path $\ga$ in $\Ga\bs X$ starting perpendicularly
from $\pi(D^-_k)$ and arriving perpendicularly to $\pi(D^+_j)$. For
every $t> 0$, we denote by $\Perp (\D^-, \D^+, t)$ the (locally
finite) set of lengths $\ell(\ga)$ at most $t$ of common
perpendiculars $\ga$ from elements of $\pi(D^-)$ to elements of
$\pi(D^+)$ (considered with multiplicities).  Let $\Perp=(\Perp (\D^-,
\D^+, N))_{N\in\NN}$.

For every $t> 0$, we denote by $\Geod(t)$ the (locally finite) set of
lengths at most $t$ of primitive closed geodesics in $\Ga\bs X$
(considered with multiplicities). If $\Ga$ is furthermore assumed to
be geometrically finite, let $\Geod=(\Geod (N))_{N\in\NN}$.

We refer to Remark (3) at the beginning of Section \ref{sec:proof} for
the use of sets with multiplicities in order to compute pair
correlations.

\bcoro\label{cor:geometric} Let $X$, $\Ga$ and $\D^\pm$ be as
above. Then the families $\Perp$ and $\Geod$ admit pair
correlation functions $g_{\Perp}$ and $g_{\Geod}$ for the scaling
function $\psi=1$ (and renormalisation to probability measures) with
$$
g_{\Perp}=g_{\Geod}:t\mapsto\frac{\delta_\Ga}{2}\;e^{-\delta_\Ga\,|t|}\;.
$$
\ecoro

\dem By \cite[Thm.~1.5]{BroParPau19}, the number of common
perpendiculars with length at most $t$ (counted with multiplicities)
is asymptotic with $\frac{\|\sigma^+_{\D^-}\|\;\|\sigma^-_{\D^+}\|}
{\delta_\Ga\;\|m_{\rm BM}\|} \; e^{\delta_\Ga\, t}$. If $\Ga$ is
furthermore assumed to be geometrically finite, by
\cite[Cor.~1.7]{PauPolSha15} and \cite[Cor.~13.5(1)]{BroParPau19}, the
number of primitive closed geodesics with length at most $t$ (counted
with multiplicities) is asymptotic with $\frac{e^{\delta_\Ga
    t}}{\delta_\Ga\,t}$ as $t\to+\infty$.  The claim follows from
Theorem \ref{theo:mainPA} with constants respectively
$(c=\frac{\|\sigma^+_{\D^-}\|\;\|\sigma^-_{\D^+}\|}
{\delta_\Ga\;\|m_{\rm BM}\|},\alpha=0,\delta=\delta_\Ga)$ and
$(c=\frac{1}{\delta_\Ga},\alpha=-1,\delta=\delta_\Ga)$.
\cqfd

\medskip 

\medskip\noindent {\bf Remarks. }  (1) Refering to \cite{PauPolSha15}
and \cite{BroParPau19} for the terminology, when $\wt F$ is a bounded
$\Ga$-invariant potential on $\Ga\bs T^1X$ which is
H\"older-continuous if $X$ is a manifold, assuming the pressure of
$\wt F$ to be positive and finite, the Gibbs measure on $\Ga\bs \G X$
for $\wt F$ to be finite and mixing for the geodesic flow, and the
outer and inner skinning measures of the families $\D^\mp$ for the
potential $\wt F$ to be finite and nonzero, then the same statement as
Corollary \ref{cor:geometric} is satisfied when $\Perp$ and $\Geod$
are endowed with weights defined by the potential as in \cite[\S
  1.2]{BroParPau19}.

\medskip
(2) The assumptions of Corollary \ref{cor:geometric} are satisfied, as
a very special case, when $X$ is a real, complex or quaternionic
hyperbolic symmetric space with finite covolume under $\Ga$, and the
images of the elements of $\D^\pm$ in $\Ga\bs X$ are points, finite
volume totally geodesic submanifolds or Margulis cusp neighbourhoods,
see \cite[Cor.~21]{ParPau17ETDS}, \cite[Theo.~3]{ParPau17MA},
\cite[Thm.~8$\cdot$1]{ParPau22MPCPS}.  For instance, if $x\in X$ and
$\D^-=\D^+=\Ga x=\{\ga x:\ga\in\Ga\}$, then $\Perp(\D^-,\D^+,t)=
\{d(x,\ga x):\ga\in\Ga\}\cap\interval 0t$, and the number of common
perpendiculars of length at most $N$ (counted with multiplicities) is
given by the growth function of the orbit $\Ga x$.

\medskip (3) In \cite{PeiTapVid20}, Peign\'e, Tapie and Vidotto
construct for all $1<\alpha<2$ examples of complete simply connected
Riemannian manifolds $X$ with pinched negative sectional curvature and
geometrically finite convergent groups $\Ga$ of isometries of $X$ such
that the growth function of the orbit of any point $x\in X$ is
asymptotic with $t\mapsto C\,t^\alpha e^{\delta_\Ga t}$ for some
$C>0$.  Theorem \ref{theo:mainPA} implies that, also in this case, the
family $\Perp$ for $\D^-=\D^+=\Ga x$ admits a pair correlation
function $g_{\Perp}$ for the scaling function $\psi=1$ (and
renormalisation to probability measures), given by $g_{\Perp}:t\mapsto
\frac{\delta_\Ga}{2}\; e^{-\delta_\Ga\,|t|}$ as in Corollary
\ref{cor:geometric}.

\medskip (4) Discrete sets with growth functions for which Theorem
\ref{theo:mainPA} can be applied to prove analogs of Corollary
\ref{cor:geometric} arise in many important dynamical systems. To name
some notable ones, Parry and Pollicott \cite{ParPol83} proved that the
number of lengths at most $t$ of closed orbits of Axiom A flows on
compact manifolds (counted with multiplicities) is asymptotic with
$t\mapsto \frac{e^{ht}}{ht}$ with $h$ the topological entropy of the
flow, and Eskin and Mirzakhani \cite{EskMir11} proved the analogous
behaviour for the lengths of closed Teichm\"uller geodesics in the
moduli space of closed Riemann surfaces of genus $g$.  Athreya,
Bufetov, Eskin and Mirzakhani \cite{AthBufEskMir12} proved the
exponential growth of orbits of the mapping class group in the
Teichm\"uller space of closed Riemann surfaces of genus $g$.

\medskip
Under additional assumptions, the asymptotic behaviour of counting functions used in the proof of Corollary \ref{cor:geometric} comes with an error term required for an application of Theorem \ref{theo:main}.

\bcoro\label{cor:geometricET} Let $X$, $\Ga$ and $\D^\pm$ be as in the
beginning of Section \ref{sec:application}.  Assume that $\Ga\bs X$ is
a compact Riemannian manifold and $m_{\rm BM}$ is exponentially mixing
under the geodesic flow for the H\"older regularity, or that $\Ga\bs
X$ is a locally symmetric space, the boundary of $D^\pm_k$ is smooth,
$m_{\rm BM}$ is finite, smooth, and exponentially mixing under the
geodesic flow for the Sobolev regularity. Assume that the strong
stable/unstable ball masses by the conditionals of $m_{\rm BM}$ are
H\"older-continuous in their radius.
 
\smallskip
Let $\psi:\NN\ra[1+\infty[$ be an at most polynomially growing scaling
function, and let $\psi':N\mapsto
\frac{\card(\Perp(\D^-,\D^+,N))^2}{\psi(N)}$ be the associated
renormalizing function. Then the family $\Perp$ has a pair correlation
function $g_{\Perp,1}:t\mapsto\frac{\delta_\Ga}{2}\;
e^{-\,\delta_\Ga\,|t|}$ if $\psi= 1$, and has Poissonian behaviour
with $g_{\Perp,\psi}= \frac{\delta_\Ga}{2}$ if
${\displaystyle\lim_{+\infty}\psi =\infty}$, with error terms as in
Theorem \ref{theo:main}.
\ecoro

\dem By \cite[Thm.~1.8 (2)]{BroParPau19}, the family $\Perp$ of common
perpendiculars has exponential growth $C\,e^{\delta_\Ga t}
(1+\bigO(e^{\kappa t}))$ for some $\kappa>0$. Thus, Theorem
\ref{theo:main} implies the claim.  \cqfd

\medskip

The geodesic flow is known to have exponential decay of H\"older
correlations for compact manifolds $M=\Ga\bs\wt M$ when $M$ is
two-dimensional by \cite{Dolgopyat98}, $M$ is $1/9$-pinched
\cite[Coro.~2.7]{GiuLivPol13}, and when $M$ is locally symmetric by
\cite{Stoyanov11}. When $X$ is a symmetric space and $\Ga$ is an
arithmetic lattice, the geodesic flow has exponential decay of Sobolev
correlations by for some $\ell\in\NN$ by
\cite[Theorem~2.4.5]{KleMar96}, with the help of \cite[Theorem
  3.1]{Clozel03} to check its spectral gap property, and of
\cite[Lemma~3.1]{KleMar99} to deal with finite cover problems. See
also \cite{MohOh15,LiTan20}.

Corollary \ref{cor:geometricET} also has generalisations when the lengths
are weighted by potentials. See, for instance, the introduction of
\cite{BroParPau19} for counting results in this generality.

{\small 
\bibliography{../biblio} 
}

\bigskip
{\small
\noindent \begin{tabular}{l} 
Department of Mathematics and Statistics, P.O. Box 35\\ 
40014 University of Jyv\"askyl\"a, FINLAND.\\
{\it e-mail: jouni.t.parkkonen@jyu.fi}
\end{tabular}
\medskip

\noindent \begin{tabular}{l}
Laboratoire de mathématique d'Orsay, UMR 8628 CNRS,\\
Bâtiment 307, Universit\'e Paris-Saclay,\\
91405 ORSAY Cedex, FRANCE\\
{\it e-mail: frederic.paulin@universite-paris-saclay.fr}
\end{tabular}}
\end{document}

%% file: macros_angl_pdf.tex
\usepackage[utf8]{inputenc}
\usepackage{amsmath,amssymb,amsthm,epsfig,bbm}
\usepackage{stmaryrd,mathabx}
\usepackage{comment}
\usepackage{color}
\usepackage[T1]{fontenc}

\usepackage{interval}

\usepackage[textsize=small]{todonotes}
\usepackage{enumitem}
\usepackage{varwidth}
\setlist{nolistsep}
\usepackage{hyperref}

\newtheorem{defi}{Definition}[section]
\newtheorem{prop}[defi]{Proposition}
\newtheorem{theo}[defi]{Theorem}
\newtheorem{theofr}[defi]{Théorème}
\newtheorem{conj}[defi]{Conjecture}
\newtheorem{lemm}[defi]{Lemma}
\newtheorem{lemmfr}[defi]{Lemme}
\newtheorem{coro}[defi]{Corollary}
\newtheorem{rema}[defi]{Remark}
\newtheorem{exem}[defi]{Example}
\newtheorem{exems}[defi]{Examples}

\newcommand{\bdefi}{\begin{defi}}
\newcommand{\edefi}{\end{defi}}
\newcommand{\bprop}{\begin{prop}}
\newcommand{\eprop}{\end{prop}}
\newcommand{\btheo}{\begin{theo}}
\newcommand{\etheo}{\end{theo}}
\newcommand{\btheofr}{\begin{theofr}}
\newcommand{\etheofr}{\end{theofr}}
\newcommand{\blemm}{\begin{lemm}}
\newcommand{\elemm}{\end{lemm}}
\newcommand{\blemmfr}{\begin{lemmfr}}
\newcommand{\elemmfr}{\end{lemmfr}}
\newcommand{\brema}{\begin{rema}}
\newcommand{\erema}{\end{rema}}
\newcommand{\bexer}{\begin{exem}}
\newcommand{\eexer}{\end{exem}}
\newcommand{\bexems}{\begin{exems}}
\newcommand{\eexems}{\end{exems}}
\newcommand{\bconj}{\begin{conj}}
\newcommand{\econj}{\end{conj}}
\newcommand{\bcoro}{\begin{coro}}
\newcommand{\ecoro}{\end{coro}}
\newcommand{\dem}{\noindent{\bf Proof. }}

\usepackage{mathrsfs}
\renewcommand\mathcal{\mathscr}

\newcommand{\D}{{\cal D}}
\newcommand{\E}{{\cal E}}
\newcommand{\F}{{\cal F}}
\newcommand{\G}{{\cal G}}

\newcommand{\N}{{\cal N}}

\newcommand{\R}{{\cal R}}

\newcommand{\maths}[1]{{\mathbb #1}}  

\newcommand{\NN}{\maths{N}}

\newcommand{\RR}{\maths{R}}

\newcommand{\ra}{\rightarrow}
\newcommand{\bs}{\backslash}

\newcommand{\wt}[1]{{\widetilde{#1}}}

\newcommand{\ga}{\gamma}
\newcommand{\Ga}{\Gamma}

\newcommand{\cqfd}{\hfill$\Box$}

\newcommand{\bigO}{\operatorname{O}}

\newcommand{\card}{{\operatorname{Card}}}

\newcommand{\Geod}{\operatorname{Geod}}

\newcommand{\Leb}{\operatorname{Leb}}

\newcommand\Perp{\operatorname{Perp}}

\newcommand{\smallo}{\operatorname{o}}

\newcommand{\Var}{\operatorname{Var}}

\newcounter{fig}

